\documentclass[reqno,11pt]{amsart}
\usepackage{amscd,amssymb,verbatim}

\numberwithin{equation}{section} \theoremstyle{plain}
\theoremstyle{plain}
\newtheorem{Thm}[subsection]{Theorem}
\newtheorem{Cor}[subsection]{Corollary}
\newtheorem{Lem}[subsection]{Lemma}
\newtheorem{prop}[subsection]{Proposition}

\theoremstyle{definition}
\newtheorem{Def}[subsection]{Definition}

\theoremstyle{remark}

\newtheorem{rem}[subsection]{Remark}

\newtheorem{exmp}{Example}

\newenvironment{thm}%
          { \begin{Thm}  }%
          { \end{Thm} }

\newenvironment{lem}%
          { \begin{Lem}    }%
          { \end{Lem} }

          { \begin{Prop}  }%
          { \end{Prop} }

          { \begin{Cor} }%
          { \end{Cor} }

          { \begin{Def} }%
          { \end{Def} }

\newcommand{\ecomp}{C_{c}(E_{s})}
\newcommand{\vcomp}{C_{c}(V)}
\newcommand{\mcomp}{C_{c}^{\infty}(M)}
\newcommand{\rncomp}{C_{c}^{\infty}}
\newcommand{\lw}{{\ell^{2}_{w}}(V)}
\newcommand{\la}{{\ell^{2}_{a}}(E_{s})}
\newcommand{\delswa}{\Delta_{\sigma}}
\newcommand{\hmax}{H_{\max}}
\newcommand{\hmin}{H_{\min}}
\newcommand{\Dom}{\operatorname{Dom}}
\newcommand{\loc}{\operatorname{loc}}

\title{A Sears-type self-adjointness result for discrete magnetic Schr\"odinger operators}

\subjclass[2000] {35J10, 39A12, 47B25}
\begin{document}
\maketitle
\begin{abstract} In the context of a weighted graph with vertex set $V$ and bounded vertex degree, we give a sufficient condition for the essential self-adjointness of the operator $\Delta_{\sigma}+W$, where $\Delta_{\sigma}$ is the  magnetic Laplacian and  $W\colon V\to\mathbb{R}$ is a function satisfying $W(x)\geq -q(x)$ for all $x\in V$, with $q\colon V\to [1,\infty)$. The condition is expressed in terms of completeness of a metric that depends on $q$ and the weights of the graph. The main result is a discrete analogue of the results of I.~Oleinik and M.~A.~Shubin in the setting of non-compact Riemannian manifolds.
\end{abstract}
\section{Introduction and the main result}\label{S:main}
\subsection{The setting}\label{SS:setting}

Let $G=(V,E)$ be an infinite graph without loops and multiple edges between vertices. By $V=V(G)$ and $E=E(G)$ we denote the set of vertices  and the set of unoriented edges of $G$ respectively. In what follows, the notation $m(x)$ indicates the degree of a vertex $x$, that is, the number of edges that meet at $x$. We assume that $G$ has bounded vertex degree: there exists a constant $N>0$ such that
\begin{equation}\label{E:assumption-bounded-deg}
m(x)\leq N,\qquad\textrm{for all }x\in V.
\end{equation}

In what follows, $x\sim y$ indicates that there is an edge that connects $x$ and $y$. We will also need a set of oriented edges
\begin{equation}\label{E: edge-or}
E_{0}:=\{[x,y],[y,x]: x,\,y\in V\textrm{ and } x\sim y\}.
\end{equation}
The notation $e=[x,y]$ indicates an oriented edge $e$ with starting vertex $o(e)=x$ and terminal vertex $t(e)=y$. The definition~(\ref{E: edge-or}) means that every unoriented edge in $E$ is represented by two oriented edges in $E_0$. Thus, there is a two-to-one map $p\colon E_0\to E$. For $e=[x,y]\in E_{0}$, we denote the corresponding reverse edge by $\widehat{e}=[y,x]$. This gives rise to an involution $e\mapsto \widehat{e}$ on $E_{0}$.

To help us write formulas in unambiguous way, we fix an orientation on each edge by specifying a subset $E_{s}$ of $E_0$ such that $E_0=E_{s}\cup \widehat{E_{s}}$ (disjoint union), where $\widehat{E_{s}}$ denotes the image of $E_{s}$ under the involution $e\mapsto \widehat{e}$. Thus, we may identify $E_{s}$ with $E$ by the map $p$.

In the sequel, we assume that $G$ is connected, that is, for any $x,\,y\in V$ there exists a path $\gamma$ joining $x$ and $y$. Here, $\gamma$ is a sequence $x_1,\,x_2,\,\dots,x_n\in V$ such that $x=x_1$, $y=x_n$, and $x_{j}\sim x_{j+1}$ for all $1\leq j\leq n-1$.


In what follows, $C(V)$ is the set of complex-valued functions on $V$, and $C(E_{s})$ is the set of functions $Y\colon E_{0}\to\mathbb{C}$ such that $Y(e)=-Y(\widehat{e})$. The notations $\vcomp$ and $\ecomp$ denote the sets of finitely supported elements of $C(V)$ and $C(E_{s})$ respectively.

In the sequel, we assume that $V$ is equipped with a weight $w\colon V\to\mathbb{R}^{+}$. By $\lw$ we denote the space of functions $f\in C(V)$ such that $\|f\|<\infty$, where $\|f\|$ is the norm corresponding to the inner product
\begin{equation}\label{E:inner-w}
(f,g):=\sum_{x\in V}w(x)f(x)\overline{g(x)}.
\end{equation}
Additionally, we assume that $E$ is equipped with a weight $a\colon E_{0}\to\mathbb{R}^{+}$ such that $a(e)=a(\widehat{e})$ for all $e\in E_{0}$. This makes $G=(G, w, a)$ a weighted graph with weights $w$ and $a$.

\subsection{Magnetic Schr\"odinger operator}\label{SS:magnetic-schro}
Let $U(1):=\{z\in\mathbb{C}\colon |z|=1\}$ and $\sigma\colon E_{0}\to U(1)$ with  $\sigma(\widehat{e})=\overline{\sigma(e)}$ for all $e\in  E_{0}$, where $\overline{z}$ denotes the complex conjugate of $z\in \mathbb{C}$.

We define the magnetic Laplacian $\delswa\colon C(V)\to C(V)$ on the graph $(G, w, a)$ by the formula
\begin{equation}\label{E:magnetic-lap}
    (\delswa u)(x)=\frac{1}{w(x)}\sum_{e\in\mathcal{O}_{x}}a(e)(u(x)-\sigma(\widehat{e})u(t(e))),
\end{equation}
where $x\in V$ and
\begin{equation}\label{E:neighborhood-o-x}
\mathcal{O}_{x}:=\{e\in E_0\colon o(e)=x\}.
\end{equation}

For the case $a\equiv 1$ and $w\equiv 1$, the definition~(\ref{E:magnetic-lap}) is the same as in~\cite{Dodziuk-Mathai-03}. For the case $\sigma\equiv 1$, see~\cite{Sunada-92} and~\cite{Torki-10}.

Let $W\colon V\to \mathbb{R}$, and consider a Schr\"odinger-type expression
\begin{equation}\label{E:magnetic-schro}
Hu:= \delswa u +Wu.
\end{equation}

Let $q\colon V\to [1,\infty)$, and assume that $W$ satisfies
\begin{equation}\label{E:minorant}
W(x)\geq -q(x), \qquad\textrm{for all }x\in V.
\end{equation}

In the sequel, we will need the notion of weighted distance on $G$. Let $w$ and $a$ be as in~(\ref{E:magnetic-lap}) and let $q$ be as in~(\ref{E:minorant}).  We define the weighted distance $d_{w,a;q}$ on $G$ as follows:
\begin{equation}\label{E:w-a-dist}
d_{w,a;q}(x,y):=\inf_{\gamma\in \Gamma_{x,y}}L_{w,a;q}(\gamma),
\end{equation}
where $\Gamma_{x,y}$ is the set of all paths $\gamma\colon x=x_1,\,x_2,\,\dots,x_n=y$ such that $x_{j}\sim x_{j+1}$ for all $1\leq j\leq n-1$, and the length $L_{w,a;q}(\gamma)$ is computed as follows:
\begin{equation}\nonumber
L_{w,a;q}(\gamma)=\sum_{j=1}^{n-1}\frac{\min\{w^{1/2}(x_j), w^{1/2}(x_{j+1})\}\cdot\min\{q^{-1/2}(x_j), q^{-1/2}(x_{j+1})\}}{\sqrt{a([x_{j},x_{j+1}])}}.
\end{equation}
In the case $q\equiv 1$, the weighted distance~(\ref{E:w-a-dist}) was defined in~\cite{vtt-10}.

We say that the metric space $(G,d_{w,a;q})$ is complete if every Cauchy sequence of vertices has a limit in $V$.

\subsection{Statement of the main result}

We now state the main result.

\begin{thm}\label{T:main-2} Assume that $(G, w, a)$ is an infinite, connected, oriented, and weighted graph. Assume that $G$ has bounded vertex degree. Assume that $W$ satisfies~(\ref{E:minorant}) and $q\colon V\to[1,\infty)$ satisfies
\begin{equation}\label{E:q-1-2-lipschitz}
|q^{-1/2}(t(e))-q^{-1/2}(o(e))|\leq C\left(\frac{\min\{w(t(e)),w(o(e))\}}{a(e)}\right)^{1/2},
\end{equation}
for all $e\in E_{s}$, where $C$ is a constant.

Additionally, assume that $(G,d_{w,a;q})$ is a complete metric space. Then, the operator $H|_{\vcomp}$ is essentially self-adjoint in $\lw$.
\end{thm}

\begin{rem} The origin of the result presented in Theorem~\ref{T:main-2} can be traced back to the paper~\cite{Sears} by D.~B.~Sears concerning the essential self-adjointness of  $(-\Delta+W)|_{\rncomp({\mathbb{R}}^{n})}$ in $L^2({\mathbb{R}}^{n})$. Here, $\Delta$ is the standard Laplacian on ${\mathbb{R}}^{n}$ and $-q\leq W\in L^{\infty}_{\loc}({\mathbb{R}}^{n})$,  where $q$ is a radially symmetric function on ${\mathbb{R}}^{n}$ satisfying properties analogous to those of Theorem 1 in the present paper (with ``completeness" replaced by the divergence of $\int_{0}^{\infty}q^{-1/2}(r)\,dr$, where $r=r(x)$ is the Euclidean distance between $x\in{\mathbb{R}}^{n}$ and $0\in{\mathbb{R}}^{n}$). We should mention that the paper~\cite{Sears} followed an idea of E.~C.~Titchmarsh~\cite{Titchmarsh49}.  More recently, I.~Oleinik~\cite{ol,Oleinik94} gave a sufficient condition for the essential self-adjointness of  $(\Delta_{M}+W)|_{\mcomp}$ in $L^2(M)$, where $\Delta_{M}$ is the scalar Laplacian on a Riemannian manifold $M$ and $-q\leq W\in L^{\infty}_{\loc}(M)$. Here, $q$ is a function on $M$ satisfying properties analogous to those of Theorem 1 in the present paper. Oleinik's proof was simplified by M.~A.~Shubin~\cite{sh1}, and the result was extended to magnetic Schr\"odinger operators in~\cite{sh2}. Theorem~\ref{T:main-2} of the present paper is a discrete analogue of the mentioned results of Oleinik and Shubin.
\end{rem}

\begin{rem} Assuming~(\ref{E:assumption-bounded-deg}), the completeness of $(G,d_{w,a;1})$, and
\begin{equation}\label{E:sa-k-sebibounded}
(Hu,u)\geq k\|u\|^2, \quad\textrm{for all }u\in\vcomp,
\end{equation}
where $k$ is a constant independent of $u$, the essential self-adjointness of $H|_{\vcomp}$  was established in~\cite[Theorem 1.3]{Milatovic-11}. If $q(x)\equiv c_0$, where $c_0$ is a constant, then the operator $H|_{\vcomp}$, with $W$ as in~(\ref{E:minorant}), satisfies~(\ref{E:sa-k-sebibounded}). However, there are operators $H$ that satisfy the hypotheses of Theorem~\ref{T:main-2} but do not satisfy~(\ref{E:sa-k-sebibounded}), as illustrated by the example below.
\end{rem}

\begin{exmp}\nonumber Consider $G=(V,E)$ with $V=\{1,2,3,\dots\}$ and $E=\{[n,n+1]\colon n\in V\}$.  Define $a([n,n+1])=1$ and $w(n)=1$, for all $n\in V$.
Let $H$ be as in~(\ref{E:magnetic-schro}) with  $\sigma([n,n+1])=1$ and $W(n)=-n^2$,  for all $n\in V$. It is is easy to see that for every $k\in\mathbb{R}$, there exists a function $u\in\vcomp$ such that the inequality~(\ref{E:sa-k-sebibounded}) is not satisfied. Thus, the operator $H$ is not semi-bounded from below, and we cannot use~\cite[Theorem 1.3]{Milatovic-11}. Turning to hypotheses of Theorem~\ref{T:main-2}, note that $W$ satisfies~(\ref{E:minorant}) with $q(n)=n^2$. It is easy to see that $q^{-1/2}=n^{-1}$ satisfies~(\ref{E:q-1-2-lipschitz}) with $C=1$. Fix $K_1\in V$, and let $K>K_1$. For $x_1=K_1$ and $x=K$, by~(\ref{E:w-a-dist}) we have
\[
d_{w,a;q}(x_1,x)=\sum_{n=K_1}^{K-1}\frac{1}{n+1}\to\infty,\quad\textrm{as }K\to\infty.
\]
Thus, the metric $d_{w,a;q}$ is complete, and by Theorem~\ref{T:main-2} the operator $H|_{\vcomp}$ is essentially self-adjoint in $\lw$.
\end{exmp}

\begin{rem}
Thanks to assumption~(\ref{E:sa-k-sebibounded}), the proof of~\cite[Theorem 1.3]{Milatovic-11} reduced to showing that if $u\in\Dom(\hmax)$, with $\hmax$ as in Section~\ref{S:preliminary} below, and $(H+\lambda)u=0$ with sufficiently large $\lambda>0$, then $u=0$. To this end, a sequence of cut-off functions was constructed and a ``summation by parts" method was used. In the absence of assumption~(\ref{E:sa-k-sebibounded}), the essential self-adjointness can be established by showing that $\hmax$ is symmetric. This requires an approach different from~\cite{Milatovic-11}: in the present paper, we consider the sum $J_{s}$ that incorporates the metric $d_{w,a;q}$ (see~(\ref{E:sum-J-s}) below) and show that $J_{s}\to 0$ as $s\to+\infty$. A key ingredient in this endeavor, not present in~\cite{Milatovic-11}, is the estimate~(\ref{E:info-dom-hmax})  for $d_{\sigma}u$, where $u\in\Dom(\hmax)$. The estimate~(\ref{E:info-dom-hmax}) is a discrete analogue of~\cite[Lemma 4.3]{sh2}.
\end{rem}

\begin{rem}
For studies of the operator~(\ref{E:magnetic-lap}) with $a\equiv 1$, $\sigma\equiv 1$, and $w\equiv m$, see, for instance,~\cite{Chung} and~\cite{Mohar-89}. For general information concerning magnetic Laplacian on graphs, see~\cite{Mathai-Yates} and~\cite{Sunada-94}. For a proof the discrete analogue of Kato's inequality, see~\cite{Dodziuk-Mathai-03}.

For the problem of self-adjoint realization of the operator~(\ref{E:magnetic-schro}) and its special cases ($a\equiv 1$, $\sigma\equiv 1$, $w\equiv 1$, and $W\equiv 0$), see, for instance, \cite{vtt-10}, \cite{vtt-10-preprint}, \cite{Golenia-Schu-10}, \cite{HKLW-preprint-11}, \cite{Jor-08-preprint}, \cite{Keller-Lenz-09}, \cite{Keller-Lenz-10}, \cite{Torki-10}, \cite{Weber-10}, and \cite{Woj-08}. We should mention that the authors of~\cite{HKLW-preprint-11} and~\cite{Keller-Lenz-09, Keller-Lenz-10} worked in the setting of discrete sets, a more general context than locally finite graphs.  For a study of the essential self-adjointness of discrete Laplace operator on forms, see~\cite{Masamune-09}.

The problem of stochastic completeness of graphs is considered in~\cite{Dodziuk-06},~\cite{Weber-10}, ~\cite{Woj-08}, and~\cite{Woj-09}. In the setting of Dirichlet forms on discrete sets, stochastic completeness is studied in~\cite{HKLW-preprint-11},~\cite{Keller-Lenz-09}, and~\cite{Keller-Lenz-10}. For another approach to stochastic completeness on discrete sets, see~\cite{Huang-11}. For a study of random walks on infinite graphs, see~\cite{Dodziuk-84},~\cite{Dodziuk-88},~\cite{Woess-00}, and references therein.

For studies of essential self-adjointness of Schr\"odinger operators in the context of non-compact Riemannian manifolds, see, for instance, \cite{br}, \cite{bms}, \cite{ga}, \cite{ol}, \cite{Oleinik94}, \cite{sh1}, \cite{sh2}, and \cite{sh01}.
\end{rem}

\section{Preliminaries}\label{S:preliminary}
In what follows, the deformed differential $d_{\sigma}\colon C(V)\to C(E_{s})$ is defined as
\begin{equation}\label{E:d-sigma}
(d_{\sigma}u)(e):=\overline{\sigma(e)}u(t(e))-u(o(e)),\qquad \textrm{for all }u\in C(V),
\end{equation}
where $\sigma$ is as in~(\ref{E:magnetic-lap}).

The deformed co-differential $\delta_{\sigma}\colon C(E_{s})\to C(V)$ is defined as
\begin{equation}\label{E:delta-sigma}
(\delta_{\sigma}Y)(x):=\frac{1}{w(x)}\sum_{\substack{e\in E_{s}\\t(e)=x}}\sigma(e)a(e)Y(e)-\frac{1}{w(x)}\sum_{\substack{e\in E_{s}\\o(e)=x}}a(e)Y(e), \end{equation}
for all $Y\in C(E_{s})$, where $\sigma$, $w$, and $a$ are as in~(\ref{E:magnetic-lap}).

In the case $\sigma\equiv 1$, the definitions~(\ref{E:d-sigma}) and~(\ref{E:delta-sigma}) give us the standard differential $d$ and standard co-differential $\delta$, respectively.

Let $\sigma$ be as in~(\ref{E:magnetic-lap}). For a function $u\in C(V)$, we define $u_{\sigma}^{\sharp}\in C(E_{s})$ by the formula
\begin{equation}\label{E:average-sigma}
    u_{\sigma}^{\sharp}(e):=\frac{\sigma(e)u(t(e))+u(o(e))}{2},\qquad\textrm{for all }e \in E_{s}.
\end{equation}

For $\sigma\equiv 1$ in~(\ref{E:average-sigma}), we define $u^{\sharp}(e):=u_{1}^{\sharp}(e)$.

In what follows, for $x\in V$, we define
\begin{equation}\label{E:neighborhood-s-x}
\mathcal{S}_{x}:=\{e\in E_s\colon o(e)=x\textrm{ or } t(e)=x\}.
\end{equation}

The proofs of the following two lemmas are straightforward computations based on the definitions of $d$, $d_{\sigma}$, $\delta$ and $\delta_{\sigma}$. For detailed proofs in the case $\sigma\equiv 1$ see~\cite[Lemma 3.1]{Masamune-09}.

\begin{lem}\label{E:leibniz-rule} For all $u\in C(V)$ and all $v\in C(V)$, the following equality holds:
\begin{equation}\label{E:derivation-sigma}
    d_{\overline{\sigma}}(uv)=(d_{\overline{\sigma}}u)v^{\sharp}+u^{\sharp}_{\sigma}dv,
\end{equation}
where $d_{\overline{\sigma}}$ is as in~(\ref{E:d-sigma})
with $\sigma(e)$ replaced by $\overline{\sigma(e)}$, $u^{\sharp}_{\sigma}$ is as in~(\ref{E:average-sigma}), and $v^{\sharp}$ is as in~(\ref{E:average-sigma}) with $\sigma\equiv 1$.
\end{lem}

\begin{lem} For all $u\in C(V)$ and all $Y\in C(E_{s})$, the following equality holds:
\begin{equation}\label{E:derivation-delta-sigma}
(\delta(u^{\sharp}_{\sigma}Y))(x)=u(x)(\delta_{\sigma}Y)(x)-\frac{1}{2w(x)}\sum_{e\in \mathcal{S}_{x}}a(e)Y(e)(d_{\overline{\sigma}}u)(e),
\end{equation}
where $d_{\overline{\sigma}}$ is as in~(\ref{E:d-sigma})
with $\sigma(e)$ replaced by $\overline{\sigma(e)}$, $u^{\sharp}_{\sigma}$ is as in~(\ref{E:average-sigma}), and $\mathcal{S}_{x}$ is as in~(\ref{E:neighborhood-s-x}).
\end{lem}

\begin{lem} Assume that $\phi\in C(V)$ is real-valued. Then
\begin{equation}\label{E:square-average}
({\phi}^{\sharp}(e))^2\leq (\phi^2)^{\sharp}(e),\qquad\textrm{for all }e\in E_{s}.
\end{equation}
\end{lem}
\noindent\textbf{Proof}
By~(\ref{E:average-sigma}) with $\sigma\equiv 1$, for all $e\in E_{s}$ we have
\[
(\phi^2)^{\sharp}(e)-({\phi}^{\sharp}(e))^2=\left(\frac{\phi(t(e))-\phi(o(e))}{2}\right)^2\geq 0,
\]
which gives~(\ref{E:square-average}).
$\hfill\square$

Let $\la$ denote the space of functions $F\in C(E_{s})$ such that $\|F\|<\infty$, where $\|F\|$ is the norm corresponding to the inner product
\[
(F,G):=\sum_{e\in E_{s}}a(e)F(e)\overline{G(e)}.
\]
It is easy to check the following equality:
\begin{equation}\label{E:adjoint-delta}
(d_{\sigma}u,Y)=(u,\delta_{\sigma}Y),\qquad\textrm{for all }u\in\lw,\,Y\in\ecomp,
\end{equation}
where $(\cdot,\cdot)$ on the left-hand side (right-hand side) denotes the inner product in $\la$ (in $\lw$).

A computation shows that the following equality holds:
\begin{equation}\label{E:laplacian-d-delta}
\delta_{\sigma}d_{\sigma}u=\Delta_{\sigma}u,\qquad\textrm{for all }u\in C(V).
\end{equation}

For the proofs of~(\ref{E:adjoint-delta}) and~(\ref{E:laplacian-d-delta}), see, for instance,~\cite[Section 3]{Milatovic-11}. The following lemma follows easily from~(\ref{E:laplacian-d-delta}) and~(\ref{E:adjoint-delta}).

\begin{lem}\label{L:symmetric} The operator $\Delta_{\sigma}|_{\vcomp}$ is symmetric in $\lw$:
\[
(\Delta_{\sigma}u,v)=(u, \Delta_{\sigma}v), \quad \textrm{for all }u,\,v\in\vcomp.
\]
\end{lem}

We now give the definitions of minimal and maximal operators associated with the expression~(\ref{E:magnetic-schro}).
We define the operator $\hmin$ by the formula
\begin{equation}\label{E:h-min}
\hmin u:=Hu,\qquad \Dom(\hmin):=\vcomp.
\end{equation}

Since $W$ is real-valued, the following lemma follows easily from Lemma~\ref{L:symmetric}.

\begin{lem}\label{L:hvmin-symm} The operator $\hmin$ is symmetric in $\lw$.
\end{lem}
We define $\hmax:=(\hmin)^{*}$, where $T^*$ denotes the adjoint of operator $T$. We also define $\mathcal{D}:=\{u\in \lw\colon Hu\in \lw\}$.

For a proof of the following lemma, see, for instance,~\cite[Lemma 3.7]{Milatovic-11}.
\begin{lem}\label{L:domain-hmax} The following hold: $\Dom(\hmax)=\mathcal{D}$ and $\hmax u=Hu$ for all $u\in\mathcal{D}$.
\end{lem}


\section{Proof of Theorem~\ref{T:main-2}}
In this section, we will adapt the technique of Shubin~\cite{sh2}.

Let $\hmin$ and $\hmax$ be as in Section~\ref{S:preliminary}. By Lemma~\ref{L:hvmin-symm} we know that $\hmin$ is symmetric. Thus, by Kato~\cite[Problem V.3.10]{Kato80}, $\hmin$ is essentially self-adjoint if and only if
\begin{equation}\label{E:h-max-symmetric}
    (\hmax u,v)=(u,\hmax v),\qquad\textrm{for all }u\,,v\in\Dom(\hmax).
\end{equation}

The following proposition provides useful information about $\Dom(\hmax)$.
\begin{prop}\label{P:info-dom-hmax} If $u\in\Dom(\hmax)$, then
\begin{align}\label{E:info-dom-hmax}
&\sum_{e\in E_{s}}\min\{q^{-1}(o(e)),q^{-1}(t(e))\}a(e)|(d_{\sigma}u)(e)|^2\nonumber\\
&\leq 2((2C^2N+1)\|u\|^2+\|Hu\|\|u\|),
\end{align}
where $H$ is as in~(\ref{E:magnetic-schro}), $N$ is as in~(\ref{E:assumption-bounded-deg}), and $C$ is as in~(\ref{E:q-1-2-lipschitz}).
\end{prop}

In the proof of Proposition~\ref{P:info-dom-hmax}, we will use a sequence  of cut-off functions. Fix a vertex $x_0\in V$, and  define
\begin{equation}\label{E:cut-off}
\chi_n(x):=\left(\left(\frac{2n-d_{w,a;1}(x_0,x)}{n}\right)\vee 0\right)\wedge 1,\qquad x\in V,\quad n\in \mathbb{Z}_{+},
\end{equation}
where  $d_{w,a;1}(x_0,x)$ is as in~(\ref{E:w-a-dist}) with $q\equiv 1$.

In the case $w\equiv 1$ and $a\equiv 1$, the sequence~(\ref{E:cut-off}) was constructed in~\cite[Proposition 3.2]{Masamune-09}. Denote
\begin{equation}\label{E:nbd-x-0}
B^{w,a}_{n}(x_0):=\{x\in V\colon d_{w,a;1}(x_0,x)\leq n\}.
\end{equation}
The sequence $\{\chi_n\}_{n\in\mathbb{Z}_{+}}$ satisfies the following properties: (i) $0\leq \chi_n(x)\leq 1$, for all $x\in V$; (ii)  $\chi_n(x)=1$ for $x\in B^{w,a}_{n}(x_0)$ and  $\chi_n(x)=0$ for $x\notin B^{w,a}_{2n}(x_0)$; (iii) for all $x\in V$, we have $\displaystyle\lim_{n\to\infty}\chi_n(x)=1$; (iv) the functions $\chi_n$ have finite support; and (v) the functions $d\chi_n$ satisfy the inequality
\begin{equation}\nonumber
|(d\chi_n)(e)|\leq \frac{d_{w,a;1}(o(e),t(e))}{n}.
\end{equation}
It is easy to see that the properties (i)--(iii) and (v) hold. By hypothesis, we know that $(G,d_{w,a;q})$ is a complete metric space and, thus, balls with respect to $d_{w,a;q}$ are finite; see, for instance,~\cite[Section 6.1]{Milatovic-11}. Let $B^{w,a;q}_{2n}(x_0)$ be as in~(\ref{E:nbd-x-0}) with $d_{w,a;1}$ replaced by $d_{w,a;q}$. Since $q\geq 1$ it follows that $B^{w,a}_{2n}(x_0)\subseteq B^{w,a;q}_{2n}(x_0)$. Thus, property (iv) is a consequence of property (ii) and the finiteness of $B^{w,a}_{2n}(x_0)$.

\noindent\textbf{Proof of Proposition~\ref{P:info-dom-hmax}}

Let $u\in\Dom(\hmax)$ and let $\phi\in\vcomp$ be a real-valued function. Define
\begin{equation}\label{E:sum-I}
I:=\left(\sum_{e\in E_{s}}a(e)|(d_{\sigma}u)(e)|^2{(\phi^2)}^{\sharp}(e)\right)^{1/2},
\end{equation}
where $f^{\sharp}(e)$ is as in~(\ref{E:average-sigma}) with $\sigma\equiv1$.

We will first show that
\begin{equation}\label{E:sum-I-show}
I^2\leq |(\phi^2Hu,u)|+(\phi^2qu,u)+2I\left(\sum_{e\in E_{s}}a(e)|(d\phi)(e)|^2|(\overline{u})_{\sigma}^{\sharp}(e)|^2\right)^{1/2},
\end{equation}
where $f_{\sigma}^{\sharp}(e)$ is as in~(\ref{E:average-sigma}), and $\overline{z}$ is the conjugate of $z\in\mathbb{C}$.

Using~(\ref{E:derivation-delta-sigma}), the equality $\Delta_{\sigma}u=Hu-Wu$, and
\begin{equation}\nonumber
(d_{\overline{\sigma}}(\phi^2\overline{u}))(e)=
\overline{(d_{\sigma}u)(e)}(\phi^2)^{\sharp}(e)+2(\overline{u})_{\sigma}^{\sharp}(e)\phi^{\sharp}(e)(d\phi)(e),
\end{equation}
we have
\begin{align}\label{E:comp-leibniz-2}
&\delta\left((\phi^2\overline{u})_{\sigma}^{\sharp}d_{\sigma} u\right)(x)=\phi^2(x)\overline{u(x)}(Hu-Wu)(x)\nonumber\\
&-\frac{1}{2w(x)}\sum_{e\in\mathcal{S}_{x}}a(e)|(d_{\sigma}u)(e)|^2(\phi^2)^{\sharp}(e)\nonumber\\
&-\frac{1}{w(x)}\sum_{e\in\mathcal{S}_{x}}a(e)(d_{\sigma}u)(e)(\overline{u})_{\sigma}^{\sharp}(e){\phi}^{\sharp}(e)(d\phi)(e).
\end{align}

Since $\phi$ has finite support, using the definition of $\delta$ it follows that
\begin{equation}\label{E:sum-stokes}
\sum_{x\in V}\left(w(x) \delta\left((\phi^2\overline{u})_{\sigma}^{\sharp}d_{\sigma} u\right)(x)\right)=0.
\end{equation}

Multiplying both sides of~(\ref{E:comp-leibniz-2}) by $w(x)$, summing over $x\in V$, and using~(\ref{E:sum-stokes}), we get
\begin{align}\label{E:comp-leibniz-3}
&\frac{1}{2}\sum_{x\in V}\sum_{e\in\mathcal{S}_{x}}a(e)|(d_{\sigma}u)(e)|^2(\phi^2)^{\sharp}(e)=(\phi^2Hu,u)-(\phi^2Wu,u)\nonumber\\
&-\sum_{x\in V}\sum_{e\in\mathcal{S}_{x}}a(e)(d_{\sigma}u)(e)(\overline{u})_{\sigma}^{\sharp}(e){\phi}^{\sharp}(e)(d\phi)(e).
\end{align}

Rewriting the double sum on the left-hand side of~(\ref{E:comp-leibniz-3}) as the sum over $E_{s}$, taking real parts on both sides of~(\ref{E:comp-leibniz-3}), and using~(\ref{E:minorant}), we have
\begin{align}\nonumber
&\sum_{e\in E_{s}}a(e)|(d_{\sigma}u)(e)|^2(\phi^2)^{\sharp}(e)=\textrm{Re }(\phi^2Hu,u)-(\phi^2Wu,u)\nonumber\\
&-\textrm{Re }\sum_{x\in V}\sum_{e\in\mathcal{S}_{x}}a(e)(d_{\sigma}u)(e)(\overline{u})_{\sigma}^{\sharp}(e){\phi}^{\sharp}(e)(d\phi)(e)\nonumber\\
&\leq |(\phi^2Hu,u)|+(\phi^2qu,u)\nonumber\\
&+2\sum_{e\in E_{s}}a(e)|(d_{\sigma}u)(e)||(\overline{u})_{\sigma}^{\sharp}(e)||{\phi}^{\sharp}(e)||(d\phi)(e)|\nonumber,
\end{align}
which, after applying Cauchy--Schwarz inequality and~(\ref{E:square-average}), gives~(\ref{E:sum-I-show}).

Let $\chi_{n}$ be as in~(\ref{E:cut-off}) and let $q$ be as in~(\ref{E:minorant}). Define
\begin{equation}\label{E:cut-off-phi-n}
\phi_{n}(x):=\chi_{n}(x)q^{-1/2}(x).
\end{equation}
By property (iv) of $\chi_n$ it follows that $\phi_n$ has finite support. By property (i) of $\chi_n$  and since $q\geq 1$, we have
\begin{equation}\label{E:cut-off-phi-n-property}
0\leq \phi_n(x)\leq q^{-1/2}(x)\leq 1,\qquad\textrm{for all }x\in V.
\end{equation}

By property (iii) of $\chi_n$  we have

\begin{equation}\label{E:cut-off-phi-n-limit}
\displaystyle\lim_{n\to\infty}\phi_n(x)=q^{-1/2}(x),\qquad\textrm{for all }x\in V.
\end{equation}

By~(\ref{E:derivation-sigma}),~(\ref{E:q-1-2-lipschitz}), properties (i) and (v) of $\chi_n$, the inequality $q\geq 1$, and~(\ref{E:w-a-dist}), we have
\begin{align}\label{E:d-phi-n}
&|(d\phi_{n})(e)|=|(d\chi_{n})(e)(q^{-1/2})^{\sharp}(e)+(\chi_{n})^{\sharp}(e)(d q^{-1/2})(e)|\nonumber\\
&\leq \left(\frac{1}{n}+C\right)\frac{\min\{w^{1/2}(o(e)),w^{1/2}(t(e))\}}{\sqrt{a(e)}},
\end{align}
where $C$ is as in~(\ref{E:q-1-2-lipschitz}).

We also have
\begin{align}\label{E:u-sharp-sigma-bound}
|(\overline{u})_{\sigma}^{\sharp}(e)|^2\leq \frac{|u(o(e))|^2+|u(t(e))|^2}{2}.
\end{align}

By~(\ref{E:d-phi-n}),~(\ref{E:u-sharp-sigma-bound}), and~(\ref{E:assumption-bounded-deg}) we get
\begin{align}\label{E:u-sharp-sigma-bound-1}
&\left(\sum_{e\in E_{s}}a(e)|(d\phi_n)(e)|^2|(\overline{u})_{\sigma}^{\sharp}(e)|^2\right)^{1/2}\nonumber\\
&\leq \frac{1}{\sqrt{2}}\left(\frac{1}{n}+C\right)\left(\sum_{e\in E_{s}}|u(o(e))|^2w(o(e))+ \sum_{e\in E_{s}}|u(t(e))|^2w(t(e))\right)^{1/2}\nonumber\\
&\leq \frac{1}{\sqrt{2}}\left(\frac{1}{n}+C\right)\left(2N\|u\|^2\right)^{1/2}=\left(\frac{1}{n}+C\right)\sqrt{N}\|u\|.
\end{align}

By~(\ref{E:sum-I-show}) with $\phi=\phi_n$,~(\ref{E:u-sharp-sigma-bound-1}), and~(\ref{E:cut-off-phi-n-property}), we obtain
\begin{equation}\label{E:quadratic-inequality}
I_n^2\leq \|Hu\|\|u\|+\|u\|^2+2I_{n}\left(\frac{1}{n}+C\right)\sqrt{N}\|u\|,
\end{equation}
for all $u\in\Dom(\hmax)$, where $I_{n}$ is as in~(\ref{E:sum-I}) with $\phi=\phi_n$.

Using the inequality $ab\leq \frac{a^2}{4}+b^2$ in the third term on the right-hand side of~(\ref{E:quadratic-inequality}) and rearranging, we obtain
\begin{equation}\label{E:quadratic-inequality-2}
I_n^2\leq 2\left(\|Hu\|\|u\|+\left(2N\left(\frac{1}{n}+C\right)^2+1\right)\|u\|^2\right).
\end{equation}

Letting $n\to\infty$ in~(\ref{E:quadratic-inequality-2}) and using~(\ref{E:cut-off-phi-n-limit}) together with Fatou's lemma, we get

\begin{equation}\label{E:quadratic-inequality-3}
\sum_{e\in E_{s}}a(e)||(d_{\sigma}u)(e)|^2(q^{-1})^{\sharp}(e)\leq 2\left(\|Hu\|\|u\|+\left(2NC^2+1\right)\|u\|^2\right).
\end{equation}
Since
\[
\min\{q^{-1}(o(e)),q^{-1}(t(e))\}\leq (q^{-1})^{\sharp}(e),\qquad\textrm{for all }e\in E_{s},
\]
the inequality~(\ref{E:info-dom-hmax}) follows directly from~(\ref{E:quadratic-inequality-3}). $\hfill\square$

\medskip

In the sequel, we will prove~(\ref{E:h-max-symmetric}). Let $d_{w,a;q}$ be as in~(\ref{E:w-a-dist}). Fix $x_0\in V$ and define
\begin{equation}\label{E:function-P}
    P(x):=d_{w,a;q}(x_0,x),\qquad x\in V.
\end{equation}

In what follows, for a function $f\colon V\to\mathbb{R}$ we define $f^{+}(x):=\max\{f(x),0\}$.

Let $u\,,v\in\Dom(\hmax)$ and let $s>0$. Define
\begin{equation}\label{E:sum-J-s}
J_{s}:=\sum_{x\in V}\left(1-\frac{P(x)}{s}\right)^{+}\left((Hu)(x)\overline{v(x)}-u(x)\overline{(Hv)(x)}\right)w(x),
\end{equation}
where $P$ is as in~(\ref{E:function-P}), $H$ is as in~(\ref{E:magnetic-schro}), and $\overline{z}$ denotes the conjugate of $z\in\mathbb{C}$.

Since $(G,d_{w,a;q})$ is a complete metric space, by~\cite[Section 6.1]{Milatovic-11} it follows that the set
\[
U_{s}:=\{x\in V\colon P(x)\leq s\}
\]
is finite. Thus, for all $s>0$, the summation in~(\ref{E:sum-J-s}) is performed over finitely many vertices.

\begin{lem}\label{L:limit-J-s} Let $J_{s}$ be as in~(\ref{E:sum-J-s}). Then
\begin{equation}\label{E:limit-J-s-symmetric}
\lim_{s\to+\infty}J_{s}=(Hu,v)-(u,Hv).
\end{equation}
\end{lem}
\noindent\textbf{Proof}
For all $x\in V$, as $s\to+\infty$, the summand in~(\ref{E:sum-J-s}) converges to
\begin{equation}\nonumber
((Hu)(x)\overline{v(x)}-u(x)\overline{(Hv)(x)})w(x).
\end{equation}
Additionally, for all $x\in V$ and $s>0$,  the summand in~(\ref{E:sum-J-s}) is estimated from above by
\begin{equation}\nonumber
|(Hu)(x)||\overline{v(x)}|w(x)+|u(x)||\overline{(Hv)(x)}|w(x).
\end{equation}
Since $u\,,v\in\Dom(\hmax)$, by Lemma~\ref{L:domain-hmax} we have $Hu\in\lw$ and $Hv\in\lw$. Hence, by Cauchy--Schwarz inequality it follows that
\begin{equation}\nonumber
\sum_{x\in V}|(Hu)(x)||\overline{v(x)}|w(x)<+\infty\quad \textrm{and}\quad \sum_{x\in V}|u(x)||\overline{(Hv)(x)}|w(x)<+\infty.
\end{equation}
Thus, by dominated convergence theorem we obtain~(\ref{E:limit-J-s-symmetric}).
$\hfill\square$

\begin{lem}\label{L:J-s-bound} Let $J_{s}$ be as in~(\ref{E:sum-J-s}) and let $N$ be as in~(\ref{E:assumption-bounded-deg}). Then
\begin{align}\label{E:J-s-bound}
&|J_s| \leq \frac{\sqrt{N}}{s}\|v\|\left(\sum_{e\in E_{s}}a(e)\min\{q^{-1}(o(e)),q^{-1}(t(e))\}|(d_{\sigma} u)(e)|^2\right)^{1/2}\nonumber\\
&+\frac{\sqrt{N}}{s}\|u\|\left(\sum_{e\in E_{s}}a(e)\min\{q^{-1}(o(e)),q^{-1}(t(e))\}|(d_{\sigma} v)(e)|^2\right)^{1/2}.
\end{align}
\end{lem}
\noindent\textbf{Proof}
Using~(\ref{E:magnetic-lap}),~(\ref{E:magnetic-schro}), and the property $\sigma(\widehat{e})=\overline{\sigma(e)}$, and recalling that $W$ is real-valued, we can rewrite~(\ref{E:sum-J-s}) as
\begin{align}\label{E:J-s-rewrite}
J_{s}=\sum_{x\in V}\sum_{e\in\mathcal{O}_{x}}\left(1-\frac{P(x)}{s}\right)^{+}a(e)
\left(\sigma(e)u(x)\overline{v(t(e))}-\sigma(\widehat{e})u(t(e))\overline{v(x)}\right).
\end{align}
An edge $e=[x,y]\in E_{0}$ occurs twice in~(\ref{E:J-s-rewrite}): once as $[x,y]$ and once as $[y,x]$.
Since $a([x,y])=a([y,x])$, it follows that the contribution of $e=[x,y]$ and $\widehat{e}=[y,x]$ together in~(\ref{E:J-s-rewrite}) is
\begin{align}\label{E:J-1-combine}
&\left(\left(1-\frac{P(x)}{s}\right)^{+}-\left(1-\frac{P(t(e))}{s}\right)^{+}\right)a(e)
\left(\sigma(e)u(x)\overline{v(t(e))}\right.\nonumber\\
&\left.-\sigma(\widehat{e})u(t(e))\overline{v(x)}\right).
\end{align}

Using~(\ref{E:J-1-combine}) and the definition of $d_{\sigma}$, we can rewrite~(\ref{E:J-s-rewrite}) as
\begin{align}\label{E:J-s-rewrite-final}
&J_{s}=\sum_{e\in E_{s}}\left(\left(1-\frac{P(o(e))}{s}\right)^{+}-\left(1-\frac{P(t(e))}{s}\right)^{+}\right)
a(e)\left(\overline{(d_{\sigma}v)(e)}u(o(e))\right.\nonumber\\
&\left.-(d_{\sigma}u)(e)\overline{v(o(e))}\right).
\end{align}
Using triangle inequality and property
\[
|f^{+}(x)-g^{+}(x)|\leq|f(x)-g(x)|,
\]
from~(\ref{E:J-s-rewrite-final}) we obtain
\begin{align}\label{E:J-s-pre-final}
&|J_{s}|\leq\frac{1}{s} \sum_{e\in E_{s}}a(e)|P(t(e))-P(o(e))|(|(d_{\sigma}v)(e)||u(o(e))|\nonumber\\
&+|(d_{\sigma}u)(e)||v(o(e))|).
\end{align}
By~(\ref{E:function-P}) and~(\ref{E:w-a-dist}) we get
\begin{align}\label{E:distances-P-d}
&|P(t(e))-P(o(e))|\leq d_{w,a;q}(t(e),o(e))\nonumber\\
&\leq \frac{w^{1/2}(o(e))\min\{q^{-1/2}(o(e)),q^{-1/2}(t(e))\}}{\sqrt{a(e)}}.
\end{align}
Combining~(\ref{E:J-s-pre-final}) and~(\ref{E:distances-P-d}), and using Cauchy--Schwarz inequality together with assumption~(\ref{E:assumption-bounded-deg}), we obtain~(\ref{E:J-s-bound}).
$\hfill\square$

\noindent\textbf{Continuation of the proof of Theorem~\ref{T:main-2}}

\medskip

Let $u\in\Dom(\hmax)$ and $v\in\Dom(\hmax)$. By Lemma~\ref{L:domain-hmax} it follows that $Hu\in\lw$ and $Hv\in\lw$. Letting $s\to+\infty$ in~(\ref{E:J-s-bound}) and using~(\ref{E:info-dom-hmax}), it follows that $J_{s}\to 0$ as $s\to+\infty$. This, together with~(\ref{E:limit-J-s-symmetric}), shows~(\ref{E:h-max-symmetric}). $\hfill\square$



\end{document}